\documentclass[12pt]{amsart}
\usepackage{times,fullpage}
\usepackage{latexsym,amscd,amssymb}
\newtheorem{thm}{Theorem}

\newtheorem{cor}{Corollary}

\theoremstyle{remark}
\newtheorem{rem}{Remark}

\theoremstyle{definition}
\newtheorem{defn}{Definition}

\newcommand{\bZ}{\mathbb Z}

\newcommand{\dvol}{{d\textrm{vol}}}

\newcommand{\Q}{\mathbb{ Q}}
\newcommand{\R}{\mathbb{ R}}

\newcommand{\Spc}{\mathrm{Spin}^c}

\begin{document}

\title{Entropies, volumes, and Einstein metrics}
\author{D.~Kotschick}
\address{Mathematisches Institut, {\smaller LMU} M\"unchen,
Theresienstr.~39, 80333 M\"unchen, Germany}
\email{dieter@member.ams.org}
\thanks{This work is part of the project {\em Asymptotic invariants of manifolds}, supported by the {\it Schwerpunktprogramm
Globale Differentialgeometrie} of the {\it Deutsche Forschungsgemeinschaft}.\\ 
Published in: 
C.~B\"ar et al.~(eds.), {\sl Global Differential Geometry}, pp.~39--54, Springer Proceedings in Mathematics {\bf 17},
Springer Verlag 2012, DOI 10.1007/978-3-642-22842-1\_ 2}

\begin{abstract}
We survey the definitions and some important properties of several asymptotic invariants 
of smooth manifolds, and discuss some open questions related to them. We prove that the 
(non-)vanishing of the minimal volume is a differentiable property, which is not invariant under 
homeomorphisms. We also formulate an obstruction to the existence of Einstein 
metrics on four-manifolds involving the volume entropy. This generalizes both the 
Gromov--Hitchin--Thorpe inequality proved in~\cite{GHT}, and Sambusetti's obstruction~\cite{S}. 
\end{abstract}

\maketitle

\section{Introduction}

In his seminal paper on bounded cohomology~\cite{gromov}, Gromov introduced 
both the simplicial volume and the minimal volume of manifolds. While the 
definition of the simplicial volume belongs to quantitative algebraic topology,
that of the minimal volume belongs to asymptotic Riemannian geometry, in
that one takes an infimum of a certain quantity over the space of all Riemannian
metrics on a fixed manifold. The simplicial volume provides a lower bound 
for the minimal volume, via the relations of these quantities to another 
asymptotic invariant, the minimal volume entropy of a manifold. This invariant 
is an asymptotic one in two ways: first one considers the asymptotics at infinity
of a quantity on a non-compact manifold, and then one takes an infimum over 
all metrics, in the same way as for the definition of the minimal volume.

In Section~\ref{s:chain} of this paper we survey several asymptotic
invariants, including the three mentioned above. The invariants we discuss are
related to each other by an interesting chain of inequalities that we explain,
with the simplicial volume at the bottom, and the minimal volume at the top.
We shall also discuss the nature of these invariants. Several of the asymptotic
Riemannian invariants actually turn out to be rather crude topological invariants,
by the recent work of Brunnbauer~\cite{br}, who completed and generalized earlier
work of Babenko~\cite{B1,B2,B3}. This is not so for the minimal volume, since
Bessi\`eres~\cite{Be} has given examples of high-dimensional manifolds that
are homeomorphic but have different minimal volumes, both of them positive. 
In Section~\ref{s:exotic} of this paper we exhibit another way in which the minimal volume 
fails to be invariant under homeomorphisms: there are pairs of homeomorphic 
four-manifolds for which the minimal volume is zero for one, and is positive 
for the other. As the vanishing of the minimal volume implies the 
vanishing of all real characteristic numbers, such examples cannot 
be simply connected. In fact, our examples are parallelizable and 
exhibit for the first time the existence of exotic smooth structures 
on parallelizable closed four-manifolds.

In Section~\ref{s:Cheeger} of this paper we discuss two further asymptotic 
invariants that can be used to bound the volume entropy, and, therefore, the minimal volume, 
from below. These invariants, based on spectral and isoperimetric quantities, are more complicated
than the invariants defined in Section~\ref{s:chain}, and it is as yet unclear whether 
Brunnbauer's results~\cite{br} about topological invariance can be adapted for them.

In Section~\ref{s:Einstein} we give an application of the discussion
of asymptotic invariants to (non-)existence of Einstein metrics 
on four-manifolds. We prove that a closed oriented four-manifold 
$X$ admitting an Einstein metric must satisfy the inequality
$$
\chi(X)\ge\frac{3}{2}\vert\sigma(X)\vert+\frac{1}{108\pi^2}\lambda(X)^{4} \ ,
$$
where $\chi(X)$ denotes the Euler characteristic, $\sigma(X)$ the 
signature, and $\lambda(X)$ the volume entropy. Using a recent rigidity result 
for the entropy due to Ledrappier and Wang~\cite{LW} we characterize 
the case of equality.

Whenever one has a positive lower bound for the volume entropy, the 
above inequality implies a particular improvement of the classical 
Hitchin--Thorpe inequality~\cite{HT}
$$
\chi(X)\ge\frac{3}{2}\vert\sigma(X)\vert \ .
$$
Two specific cases of lower bounds for the entropy which we discuss in 
Section~\ref{s:Einstein} are the following. First, from Section~\ref{s:chain}, 
the spherical volume and the simplicial volume give lower bounds for 
the entropy. In the case of the simplicial volume this leads to
$$
\chi(X)\geq\frac{3}{2}\vert\sigma(X)\vert+\frac{1}{162\pi^2}\vert\vert 
X\vert\vert \ , 
$$
which improves on our result in~\cite{GHT}. Second, using the lower 
bound for the entropy arising from the existence of maps to locally 
symmetric spaces of rank one proved by Besson--Courtois--Gallot~\cite{BCG}, 
we deduce the results of Sambusetti~\cite{S}. Thus, our obstruction to 
the existence of Einstein metrics  
subsumes all the known obstructions which are homotopy invariant. 
See~\cite{K,monopoles} for obstructions which are not homotopy invariant.

\section{Entropies and volumes}\label{s:chain}

Let $M$ be a connected closed oriented manifold of dimension $n$. In this 
section we discuss the following chain of inequalities between its topological 
invariants:
\begin{equation}\label{e:main}
    \frac{n^{n/2}}{n!} \vert\vert M \vert\vert \leq 2^{n}n^{n/2} T(M) \leq 
    \lambda(M)^{n}\leq h(M)^{n}\leq (n-1)^{n}\textrm{MinVol}(M) \ .
    \end{equation}
This chain appears in~\cite{PP}, without $T(M)$, and with an unspecified 
constant in front of $\vert\vert M \vert\vert$. The inequalities involving 
$T(M)$ are from~\cite{BCG0}, where $T(M)$ was defined. 

We now explain the terms in~(\ref{e:main}), and the sources for the 
different inequalities.

\subsection{Simplicial volume}
The simplicial volume was introduced by Gromov~\cite{gromov}.
Let $c=\Sigma_i r_i \sigma_i$ be a chain with real coefficients $r_i$, 
where $\sigma_i\colon\Delta^k\rightarrow X$ are singular $k$--simplices in $M$.
One defines the norm of $c$ to be
$$
\vert\vert c\vert\vert =\Sigma_i\vert r_i\vert \ .
$$
If $\alpha\in H_k(M,\R )$, set
$$
\vert\vert\alpha\vert\vert = inf \{\vert\vert c\vert\vert \ | \ 
c \ \textrm{a cycle representing} \ \alpha\} \ ,
$$
where the infimum is taken over all cycles representing $\alpha$.
This is a seminorm on the homology. 

The simplicial volume is defined to be the norm of the fundamental class, 
and is denoted
$$
\vert\vert M\vert\vert = \vert\vert [M]\vert\vert \ .
$$

Gromov~\cite{gromov} has shown that the simplicial volume of $M$ is 
determined by the image of its fundamental class under the classifying 
map $f\colon M\rightarrow B\pi_{1}(M)$ of the universal covering. 

This motivates the following definition.

\begin{defn}\label{d:hominv}
An invariant $I$ of closed oriented manifolds $M$ is said to be homologically invariant
if it depends only on $f_*[M]\in H_*(B\pi_1(M))$.
\end{defn}

Homological invariance implies that $I$ is homotopy-invariant, and that it 
is a bordism invariant of $[M,f]\in\Omega_{n}(B\pi_{1}(M))$.

\subsection{Spherical volume}
The spherical volume was introduced by Besson--Courtois--Gallot~\cite{BCG0}.

Let $\tilde M$ be the universal covering of $M$. Fix a positive 
$\pi_{1}(M)$-invariant measure $\tilde\mu$ on $\tilde M$ that is 
absolutely continuous with respect to Lebesgue measure. Denote by 
$(S^{\infty},can)$ the unit sphere in the Hilbert space 
$L^{2}(\tilde M , \tilde\mu)$ endowed with the metric induced by the 
scalar product. For every $\pi_{1}(M)$-equivariant immersion 
$\Phi\colon\tilde M\rightarrow S^{\infty}$ we have an induced metric 
$\Phi^{*}(can)$ which descends to $M$, and we write $vol(\Phi)$ for 
$Vol(M, \Phi^{*}(can))$.

The spherical volume of $M$ is defined to be 
$$
T(M) = inf\{ vol(\Phi) \ 
\vert \ \Phi \ \textrm{a} \ \pi_{1}(M)-\textrm{equivariant immersion} 
\ \tilde M\rightarrow S^{\infty}\} \ .
$$
This does not depend on the choice of measure $\tilde\mu$, see~\cite{BCG0}.

Brunnbauer~\cite{br} has proved that the spherical volume is homologically 
invariant in the sense of Definition~\ref{d:hominv}. In particular it is homotopy-invariant.

Besson, Courtois and Gallot proved the first inequality in~(\ref{e:main}) as 
Theorem 3.16 in~\cite{BCG0}.

\subsection{Volume entropy}
For a Riemannian metric $g$ on $M$ consider the lift $\tilde g$ to the 
universal covering $\tilde M$. For an arbitrary basepoint $p\in\tilde M$ 
consider the limit 
$$
\lambda (M,g) = \lim_{R\to\infty}\frac{log Vol(B(p,R))}{R} \ ,
$$
where $B(p,R)$ is the ball of radius $R$ around $p$ in $\tilde M$ with respect 
to $\tilde g$, and the volume is taken with respect to $\tilde g$ as well. After 
earlier work by Efremovich, Shvarts, Milnor~\cite{Mi} and others, Manning~\cite{M} 
showed that the limit exists and is independent of $p$. It follows from~\cite{Mi} 
that $\lambda (M,g)>0$ if and only if $\pi_{1}(M)$ has exponential growth. 

We call $\lambda (M,g)$ the volume entropy of the metric $g$, and define the 
volume entropy of $M$ to be 
$$
\lambda (M) = inf\{\lambda (M,g)\ \vert \  g\in Met(M) \ 
\textrm{with} \ Vol(M,g)=1 \} \ .
$$
This sometimes vanishes even when $\lambda (M,g)>0$ for every $g$.
The normalization of the total volume is necessary because of the 
scaling properties of $\lambda (M,g)$.

Any metric on $M$ can be scaled so that $Ric_g \geq -\frac{1}{n-1}g$.
By the Bishop volume comparison theorem, cf.~\cite{BC}, this implies 
$\lambda (M,g)\leq 1$. Setting
$$
\hat g = \frac{1}{Vol(M,g)^{2/n}}g
$$
we have $Vol(M,\hat g)=1$ and $\lambda (M,\hat g) = Vol(M,g)^{1/n}\lambda 
(M,g)\leq Vol(M,g)^{1/n}$, implying
\begin{equation}\label{e:Bishop}
    \lambda (M) = \inf\{\lambda (M,\hat g)\}
    \leq \inf\{ 
    Vol(M,g)^{1/n} \ \vert \ g\in Met(M) \ \textrm{s.t.} 
    \ Ric_g \geq -\frac{1}{n-1}g \}.
    \end{equation}
This is not an equality because the Bishop estimate is not sharp, 
except for metrics of constant sectional curvature.

Babenko has shown that the volume entropy $\lambda (M)$ is homotopy 
invariant~\cite{B1}, and is also an invariant of the bordism class 
$[M,f]\in\Omega_{n}(B\pi_{1}(M))$, cf.~\cite{B3}. These results were 
sharpened by Brunnbauer~\cite{br}, who proved the homological invariance of $\lambda (M)$.
Although these results show that $\lambda (M)$ is a rather crude invariant
which is understandable in many situations, there are still many open
questions about this invariant. For example, its behavior under taking 
finite coverings is not understood. We do not even know whether the vanishing
of $\lambda$ on a finite covering of $M$ implies the vanishing on $M$ itself.
This question came up in~\cite{KKM}.

The second inequality in~(\ref{e:main}) follows from Theorem~3.8 of
Besson--Courtois--Gallot~\cite{BCG0}.

\subsection{Topological entropy}
For a Riemannian metric $g$ on $M$ consider the topological entropy 
$h(M,g)$ of its geodesic flow as a dynamical system on the unit sphere 
bundle, cf.~\cite{M}. The topological entropy of $M$ is defined to be
$$
h (M) = inf\{h (M,g)\ \vert \  g\in Met(M) \ 
\textrm{with} \ Vol(M,g)=1 \} \ .
$$
Here again the normalization of the total volume is necessary because of  
scaling properties.

It seems to be unkown what exactly the minimal topological entropy depends on,
e.g. the homotopy type, the homeomorphism type, or the diffeomorphism type.
It is not clear whether this is a subtle invariant, or a crude one like the minimal volume 
entropy and the spherical volume.

Manning~\cite{M} proved $\lambda (M,g)\leq h(M,g)$ for all closed $M$.
Taking the infimum over all metrics with normalized volume yields the 
third inequality in~(\ref{e:main}).

\subsection{Minimal volume}
The minimal volume was introduced by Gromov~\cite{gromov}. It is 
defined by 
$$
\textrm{MinVol}(M)=inf\{ Vol(M,g) \ \vert \ g\in Met(M) \ 
\textrm{with} \ \vert K_{g}\vert\leq 1 \} \ , 
$$
where $K_{g}$ denotes the sectional curvature of $g$. 

It is known that the minimal volume is a very sensitive invariant of 
$M$, which depends on the smooth structure in an essential way. 
Bessi\`eres~\cite{Be} has given examples of pairs of high-dimensional 
manifolds which are homeomorphic, but have different positive minimal 
volumes. In Section~\ref{s:exotic} we shall use Seiberg--Witten theory 
to show that the vanishing of the minimal volume is not invariant under 
homeomorphisms. This represents a more dramatic failure of topological 
invariance than Bessi\`eres' result, since it shows that there are homeomorphic
manifolds such that one collapses with bounded sectional curvature, and 
the other one does not. 

Manning~\cite{M2} proved that for a closed Riemannian manifold with sectional 
curvature bounded by $\vert K_{g}\vert\leq k$ the topological entropy 
is bounded by
$$
h(M,g)\leq (n-1)\sqrt{k} \ .
$$
By rescaling, this implies the last of the inequalities in~(\ref{e:main}), 
cf.~\cite{P} p.~129.

\section{Isoperimetric constants and minimal eigenvalues}\label{s:Cheeger}

One of Gromov's original motivations for introducing the simplicial volume was to obtain lower
bounds for the minimal volume. However, the lower bounds given by the simplicial volume via~(\ref{e:main}) are
usually rather weak. In particular, there seems to be no known example of a manifold with 
$T(M)>0$ for which the inequality
\begin{equation}\label{e:T}
 \frac{n^{n/2}}{n!} \vert\vert M \vert\vert \leq 2^{n}n^{n/2} T(M)
\end{equation}
is sharp. All the other inequalities in~(\ref{e:main}) are 
sharp for surfaces of genus $g\geq 2$, where 
$$
4(g-1)=\vert\vert\Sigma_{g}\vert\vert <8T(\Sigma_{g}) =
\lambda(\Sigma_{g})^{2}=h(\Sigma_{g})^{2}=\textrm{MinVol}(\Sigma_{g})=4\pi(g-1) \ .
$$
One approach to understanding and quantifying the failure of~(\ref{e:T}) to be sharp
is to consider intermediate invariants which may interpolate between the simplicial and 
the spherical volumes. 
In fact, a candidate for such an intermediate invariant, based on the minimal eigenvalue of the Laplacian
on the universal covering, occurs in the work of Besson--Courtois--Gallot~\cite{BCG0}.
We now elaborate on this to discuss the following alternative to~(\ref{e:main}):
\begin{equation}\label{e:mainalt}
I(M)^n\leq 2^n\Lambda_0(M)^{n/2}   \leq 2^{n}n^{n/2} T(M) \leq  \lambda(M)^{n}\leq h(M)^{n} \ ,
    \end{equation}
where we can also continue on the right with the same minimal volume term as  in~(\ref{e:main}).
Only the first two terms and the first two inequalities on the left need any explanation, as the rest
has been explained already in the previous section.

\subsection{Minimal eigenvalue}
Given a closed Riemannian manifold $(M,g)$, we consider the Riemannian universal cover
$(\tilde M,\tilde g)$, and define
$$
\lambda_0(\tilde M,\tilde g) = 
\inf_{f\in C^{\infty}_{0}}\frac{\int_{\tilde M}f\cdot \Delta f \ dvol_{\tilde g}}{\int_{\tilde M}f^2 dvol_{\tilde g}} \ ,
$$
where $C^{\infty}_{0}$ denotes the smooth compactly supported functions on $\tilde M$.
Extending the Laplacian to $L^2$-functions, $\lambda_0(\tilde M,\tilde g)$ is the greatest lower 
bound for its spectrum. It is tempting now to define the minimal eigenvalue of $M$ as
$$
\inf\{\lambda_0 (\tilde M,\tilde g)\ \vert \  g\in Met(M) \ 
\textrm{with} \ Vol(M,g)=1 \} \ ,
$$
which would be the naive generalization of the definition of the minimal volume entropy. However, this 
definition always gives zero, cf.~\cite{Sch}. The correct definition of an invariant of $M$
derived from the minimal eigenvalue is the following. First, we define an invariant of a conformal
class $[g]$ by setting
$$
\Lambda_0(M,[g]) =\sup\{\lambda_0 (\tilde M,\tilde g')\ \vert \  g'\in [g] \ 
\textrm{with} \ Vol(M,g')=1 \} \ ,
$$
and then define
$$
\Lambda_0(M)=\inf\{\Lambda_0 (M,[g])\ \vert \  [g]\in Conf(M)  \} 
$$
as an infimum over conformal classes. 

This is a meaningful invariant, which for surfaces of genus $g\geq 2$ can be shown to equal $\pi (g-1)$, 
cf.~\cite{Sch}. The work of Besson--Courtois--Gallot~\cite{BCG0} implies the second inequality in~(\ref{e:mainalt})
in all dimensions $n$. It is sharp for surfaces.

\subsection{Isoperimetric constant}
Cheeger's isoperimetric constant of $(M,g)$ is defined to be
$$
i(\tilde M,\tilde g) = \inf_N\frac{Vol_{n-1}(\partial N)}{Vol_n(N)} \ ,
$$
where $N\subset \tilde M$ ranges over all compact connected subsets with smooth boundary, say. 
Again one would naively define the isoperimetric constant of $M$ to be
$$
\inf\{i (\tilde M,\tilde g)\ \vert \  g\in Met(M) \ 
\textrm{with} \ Vol(M,g)=1 \} \ ,
$$
but this always gives zero, cf.~\cite{Sch}. The definition has to  be modified in the same way as for the minimal 
eigenvalue. First we define an invariant of a  conformal class $[g]$ by setting
$$
I(M,[g]) = \sup \{i (\tilde M,\tilde g')\ \vert \  g'\in [g] \ 
\textrm{with} \ Vol(M,g')=1 \} \ ,
$$
and then we take the infimum over conformal classes:
$$
I(M) = \inf \{ I(M,[g])  \ \vert \ [g]\in Conf(M)\} \ .
$$

Cheeger's inequality $i (\tilde M,\tilde g)^2\leq 4\lambda_0(M)$, see~\cite{Cheeger}, implies the 
first inequality in~(\ref{e:mainalt}).

\bigskip

For both the minimal eigenvalue $\Lambda_0(M)$ and for the isoperimetric constant $I(M)$ it 
is not immediately clear what they depend on. One could try to use Brunnbauer's axiomatic approach~\cite{br}
to prove that these are homological invariants of manifolds. However, because of the  complicated  
definitions, using a supremum within each conformal class before taking the infimum over conformal classes, it 
is hard to verify the required axioms for these invariants. Because of this difficulty we do not yet know 
how subtle these invariants are.

The chain~(\ref{e:mainalt}) shows that one can use the minimal eigenvalue and the isoperimetric constant instead of the simplicial
volume to obtain lower bounds for the minimal volume entropy and the minimal volume. 

A very important outstanding problem is to decide whether there is an upper bound for the simplicial 
volume in terms of the isoperimetric constant. More precisely, one would like to know whether the inequality 
$$
 \frac{n^{n/2}}{n!} \vert\vert M \vert\vert \leq I(M)^n
$$
holds. If this were the case, then we could insert~(\ref{e:mainalt}) into~(\ref{e:main}), and we could hope to
measure the gap between the simplicial volume and the spherical volume. Potentially this could lead
to an improvement of~(\ref{e:main}), by improving the constant in front of the simplicial volume to close 
the gap between the simplicial and spherical volumes in~(\ref{e:main}).

\section{Einstein metrics on four-manifolds}\label{s:Einstein}

We now prove the following constraint on the topology of four-manifolds 
admitting Einstein metrics.
\begin{thm}\label{t:GHT}
Let $X$ be a closed oriented Einstein $4$--manifold. Then
\begin{equation}\label{eq:GHT}
\chi(X)\ge\frac{3}{2}\vert\sigma(X)\vert+\frac{1}{108\pi^2}\lambda(X)^{4} \ ,
\end{equation}
where $\chi(X)$ denotes the Euler characteristic, $\sigma(X)$ the 
signature, and $\lambda(X)$ the volume entropy. 

Equality in~(\ref{eq:GHT}) occurs if and only if every Einstein metric 
on $X$ is flat, is non-flat locally Calabi--Yau, or is of constant negative 
sectional curvature.
\end{thm}

\begin{rem}\label{r:B}
Note that the right-hand side of~(\ref{eq:GHT}) is an invariant of the 
cobordism class $[X,f]\in\Omega_{4}(B\pi_{1}(X))$, where 
$f\colon X\rightarrow B\pi_{1}(X)$ is the classifying map of the universal 
covering. 
For the signature this is due to Thom, and for the volume entropy it was 
proved by Babenko\footnote{His paper assumes that manifolds are of 
dimension $\geq 5$, but that is not important here; cf.~\cite{br}.} in~\cite{B3}. 
Rationally, one has the well-known isomorphism
$$
\Omega_{4}(B\pi_{1}(X))\otimes\Q = H_{4}(B\pi_{1}(X);\Omega_{0}(\star)\otimes\Q )
\oplus H_{0}(B\pi_{1}(X);\Omega_{4}(\star)\otimes\Q ) \ .
$$
The second summand is responsible for the signature term (=first Pontryagin 
number) and the first one for the entropy term in~(\ref{eq:GHT}). Indeed, Brunnbauer's 
homological invariance result for the entropy~\cite{br} shows that $\lambda(X)$
depends only on $f_{*}[X]\in H_{4}(B\pi_{1}(X))$.
\end{rem}

\begin{proof}[Proof of Theorem~\ref{t:GHT}]
By the Gauss--Bonnet theorem the Euler characteristic 
of a closed oriented Riemannian $4$--manifold $(X,g)$ is
$$
\chi(X)=\frac{1}{8\pi^2}\int_X \frac{1}{24}s_g^2 +|W|^2-|Ric_0|^2
 \dvol_{g} \ ,
$$
where $s_g$ is the scalar curvature, $W$ is the Weyl tensor,
and $Ric_0$ is the trace-less Ricci tensor of $g$.
For an Einstein metric this reduces to
\begin{equation}\label{eq:GB}
\chi(X)=\frac{1}{8\pi^2}\int_X \frac{1}{24}s_g^2+
\vert W\vert^2 \dvol_g \ .
\end{equation}
Thus $\chi(X)\geq 0$, with equality if and only if every Einstein 
metric on $X$ is flat.

If the scalar curvature of an Einstein metric is positive, then 
$\pi_1(X)$ is finite by Myers's theorem and the volume entropy 
vanishes, so that~(\ref{eq:GHT}) reduces to the Hitchin--Thorpe 
inequality obtained by comparing~(\ref{eq:GB}) with the Chern--Weil 
formula
\begin{equation}\label{eq:CW}
\sigma (X)=\frac{1}{12\pi^2}\int_X |W_+|^2-|W_-|^2 \dvol_{g} \ .
\end{equation}

If the scalar curvature is zero, then the Cheeger--Gromoll splitting 
theorem implies that either the fundamental group is finite, or the 
Euler characteristic vanishes. In the latter case the Einstein metric 
is flat and the volume growth is only polynomial, so that in both cases 
the volume entropy vanishes and we are done as before.

If the scalar curvature is negative, we scale the metric so that
$Ric_g = -\frac{1}{3}g$. 

Using the Chern--Weil formula~(\ref{eq:CW}) the second term  
in~(\ref{eq:GB}) is
\begin{equation}\label{eq:second}
\frac{1}{8\pi^2}\int_X\vert W\vert^2 \dvol_g \geq\frac{1}{8\pi^2}\vert
\int_X\vert W_+\vert^2-\vert W_-\vert^2 \dvol_g \vert =\frac{3}{2}
\vert\sigma (X)\vert \ .
\end{equation}

The first term in~(\ref{eq:GB}) is
$$
\frac{1}{8\pi^2}\frac{1}{24}\left(-\frac{4}{3}\right)^2 Vol(X,g) = 
\frac{1}{108\pi^2}Vol(X,g) \ . 
$$
Now using the Bishop estimate to bound the volume entropy $\lambda(X,g)$ from
 above by $1$ (see~(\ref{e:Bishop})) we find
 \begin{equation}\label{eq:first}
\frac{1}{8\pi^2}\frac{1}{24}\left(-\frac{4}{3}\right)^2 Vol(X,g)
\ge \frac{1}{108\pi^2} Vol(X,g) \lambda(X,g)^4
\ge \frac{1}{108\pi^2}\lambda(X)^{4} \ .
\end{equation}
This completes the proof of~(\ref{eq:GHT}). Equality cannot hold in the 
case of positive scalar curvature because in this case we threw away the 
scalar curvature term in~(\ref{eq:GB}). In the case of zero scalar 
curvature the entropy vanishes and the discussion of equality reduces 
to the corresponding discussion for the Hitchin--Thorpe inequality, see 
Hitchin~\cite{HT}. The conclusion is that every Einstein metric is flat 
or (up to choosing the orientation suitably) locally Calabi--Yau but 
non-flat. For negative scalar curvature equality in~(\ref{eq:GHT})
implies equality in~(\ref{eq:first}), so that $\lambda(X)>0$ and every 
Einstein metric on $X$ has to be entropy-minimizing. Moreover, the Einstein
metric must have entropy $\lambda(X,g)=1$, and so by the recent 
rigidity theorem of Ledrappier and Wang~\cite{LW} it has 
constant negative sectional curvature\footnote{Our scaling normalization for
the lower Ricci curvature bound is 
different from the one used in~\cite{LW}, so entropy $=1$ in our case corresponds
to entropy $=n-1$ in~\cite{LW}.}. Thus these 
are the only candidates for equality. They do indeed give equality 
because they are conformally flat and therefore~(\ref{eq:second}) is 
an equality for them, and they are entropy-minimizing by 
the celebrated result of Besson--Courtois--Gallot~\cite{BCG}.
\end{proof}

\begin{rem}
    The three cases giving rise to equality in~(\ref{eq:GHT}) correspond 
    to the vanishing of all three terms in~(\ref{eq:GHT}) for flat manifolds, 
    to the vanishing of $\lambda(X)$ only for Calabi--Yaus, and to the 
    vanishing of $\sigma(X)$ only in the hyperbolic case. After earlier, 
    unpublished, work of Calabi, Charlap--Sah, and Levine, the closed 
    orientable flat four-manifolds were classified by Hillman~\cite{Hil} and 
    by Wagner~\cite{wagner}, who showed that there are $27$ distinct ones. By 
    Bieberbach's theorems, all these manifolds are finite quotients of $T^{4}$. 
    In the locally Calabi--Yau case, Hitchin~\cite{HT} showed that the 
    only possible manifolds are the $K3$ surface and its quotients by 
    certain free actions of $\bZ_{2}$ and of $\bZ_{2}\times\bZ_{2}$. 
    In the hyperbolic case there are of course infinitely many manifolds.
    \end{rem}

Combining~(\ref{e:mainalt}) with Theorem~\ref{t:GHT} we obtain:
\begin{cor}
Let $X$ be a closed oriented Einstein $4$--manifold. Then we have the 
following lower bounds for the Euler characteristic of $X$:
\begin{equation}\label{eq:GHTT}
\chi(X)\ge \frac{3}{2}\vert\sigma(X)\vert+\frac{64}{27\pi^2} T(X) \ ,
\end{equation}    
\begin{equation}\label{eq:GHTL}
\chi(X)\ge \frac{3}{2}\vert\sigma(X)\vert+\frac{4}{27\pi^2} \Lambda_0(X)^2 \ ,
\end{equation}    
\begin{equation}\label{eq:GHTI}
\chi(X)\ge \frac{3}{2}\vert\sigma(X)\vert+\frac{1}{108\pi^2} I(X)^4 \ .
\end{equation}    
In all cases, equality occurs if and only if every Einstein metric 
on $X$ is flat, is non-flat locally Calabi--Yau, or is of constant negative 
sectional curvature.
    \end{cor}
The point of the last statement is that equality in one of these estimates implies equality in (\ref{eq:GHT}).

Combining~(\ref{e:main}) with Theorem~\ref{t:GHT} we obtain the following 
improvement of the Gromov--Hitchin--Thorpe inequality proved in~\cite{GHT}:
\begin{cor}
Let $X$ be a closed oriented Einstein $4$--manifold. Then
\begin{equation}\label{eq:GHT2}
\chi(X)\geq
\frac{3}{2}\vert\sigma(X)\vert+\frac{1}{162\pi^2}\vert\vert 
X\vert\vert \ .
\end{equation}    
    \end{cor}
Note that this is not sharp in any interesting cases with non-vanishing simplicial volume term, because
the lower bound for the entropy in terms of the simplicial volume is not optimal.

Theorem~\ref{t:GHT} together with the work of Besson--Courtois--Gallot~\cite{BCG} 
implies the following:
\begin{cor}[Sambusetti~\cite{S}]
    Let $X$ be a closed oriented Einstein $4$-manifold. 
    
    If $X$ admits a map of non-zero degree $d$ to a real hyperbolic manifold $Y$, then 
    \begin{equation}\label{eq:real}
    \chi(X)\ge\frac{3}{2}\vert\sigma(X)\vert+\vert d\vert\cdot\chi(Y) \ ,
    \end{equation}
    with equality only if $X$ is also real hyperbolic.
    
    If $X$ admits a map of non-zero degree $d$ to a complex hyperbolic manifold $Y$, then 
    \begin{equation}\label{eq:cx}
    \chi(X)>\frac{3}{2}\vert\sigma(X)\vert+\vert d\vert\cdot \frac{32}{81}\chi(Y) \ .
    \end{equation}
    \end{cor}
    \begin{proof}
	To prove~(\ref{eq:real}) we combine~(\ref{eq:GHT}) with the 
	inequality
	$$
	\lambda(X)^{4}\geq d\cdot\lambda (Y)^{4}=d\cdot 3^{4}Vol(Y,g_{0})
	$$ 
	from~\cite{BCG}. Here $g_{0}$ is the hyperbolic metric on $Y$ 
	normalized so that $K=-1$. With this normalization the Gauss--Bonnet 
	formula~(\ref{eq:GB}) gives $Vol(Y,g_{0})=\frac{4\pi^{2}}{3}\chi(Y)$ 
	so that $3^{4}Vol(Y,g_{0})=108\pi^{2}\chi(Y)$.
	
	In the case of equality $X$ must itself be hyperbolic, by the 
	equality case of Theorem~\ref{t:GHT}. 
	
	For maps to complex hyperbolic manifolds one uses the corresponding statement 
	from~\cite{BCG}. Again $\lambda(X)^{4}\geq d\cdot\lambda (Y)^{4}$ and the 
	hyperbolic metric on $Y$ is entropy-minimizing. Thus one only has to check 
	the proportionality factors between the fourth power of the entropy and the 
	Euler characteristic, and the weak form of~(\ref{eq:cx}) follows. By 
	Theorem~\ref{t:GHT} equality cannot occur in this case. 
	\end{proof}
	
It would be interesting to know whether~(\ref{eq:GHT}) remains true if we replace 
the volume entropy $\lambda$ by the topological entropy $h$, cf.~(\ref{e:main}). 
The results of Paternain and Petean~\cite{PPMRL,PP} are very suggestive in this regard. 
While the volume entropy can only be positive for manifolds with fundamental groups of 
exponential growth, the topological entropy may be positive even for simply connected 
manifolds. Unlike the volume entropy, the topological entropy is not known to be 
homotopy invariant. Notice that we can definitely not replace $\lambda^{4}$ by a 
positive multiple of $\textrm{MinVol}$, because the $K3$ surface satisfies $\chi(K3)=
\frac{3}{2}\vert\sigma(K3)\vert$, and has positive $\textrm{MinVol}$ as its Euler 
characteristic and signature are non-zero. We will show in the next section that 
the non-vanishing of the minimal volume depends in an essential way on the smooth 
structure. That the existence of an Einstein metric depends on the smooth structure was 
first shown in~\cite{K}. Recently Brunnbauer, Ishida and Su\'arez-Serrato~\cite{BIS}
have constructed some interesting examples showing that the smooth obstructions
to the existence of Einstein metrics are completely independent of those provided by 
Theorem~\ref{t:GHT}.

To end this section, we briefly discuss the history of Theorem~\ref{t:GHT}. The Hitchin--Thorpe
inequality $\chi(X)\geq\frac{3}{2}\vert\sigma (X)\vert$ was, for a long time,
the only known obstruction to the existence of Einstein metrics on four-manifolds. This changed when
Gromov~\cite{gromov} gave a lower bound for the Euler number by a multiple of the simplicial
volume. For many years after that, the Hitchin--Thorpe inequality and the ``Gromov obstruction''
were treated as separate, unrelated obstructions; see for example the discussion in~\cite{besse,BCG,S}.
In 1997, trying to understand Gromov's argument, I found that the signature term of the Hitchin--Thorpe
inequality and the simplicial volume term of Gromov's inequality can actually be combined, to 
obtain an inequality like (\ref{eq:GHT2}); see~\cite{GHT}. That inequality still has the flaw that it is 
not sharp in any interesting cases, because the simplicial volume term is too 
weak. Thinking about the chain~(\ref{e:main}), I found the sharp Theorem~\ref{t:GHT} in 2004, 
wrote it down in~\cite{first}, and also explained that it subsumes all the known homotopy-invariant 
obstructions to Einstein metrics. That paper has remained unpublished, because various referees claimed 
that Theorem~\ref{t:GHT} was not interesting, or that it was the same as the result of~\cite{GHT}, or 
that it was known to the authors of~\cite{BCG,S}. What the referees did not
notice, and I myself only noticed recently, is that the discussion of the limiting case in~\cite{first} was 
actually incomplete. The only way I know how to characterize the limiting case, is through the recent 
rigidity result of Ledrappier and Wang~\cite{LW}, as used in the proof of Theorem~\ref{t:GHT} given
above. In particular there seems to be no way of obtaining the desired conclusion from~\cite{BCG,S}.

\section{Minimal volumes and smooth structures}\label{s:exotic}

In this section we show that vanishing of the minimal volume is a property of the 
smooth structure, which is not invariant under homeomorphisms. The proof below 
actually shows that one can change the smooth structure of a manifold with a 
smooth free circle action so that for the new smooth structure any smooth
circle action must have fixed points.
\begin{thm}\label{t:2}
    For every $k\geq 0$ the manifold $X_{k}=k(S^{2}\times S^{2})\# (1+k)(S^{1}\times 
	S^{3})$ with its standard smooth structure has zero minimal volume.
	
	If $k$ is odd and large enough, then there are infinitely many pairwise 
	non-diffeomorphic smooth manifolds $Y_{k}$ homeomorphic to $X_{k}$, all of 
	which have strictly positive minimal volume.
    \end{thm}
    \begin{proof}
	Note that $X_{0}=S^{1}\times S^{3}$ has obvious free circle actions, and 
	therefore collapses with bounded sectional curvature. To see that all 
	$X_{k}$ have vanishing minimal volume it suffices to construct 
	fixed-point-free circle actions on them.
	
	The product $S^{2}\times S^{2}$ has a diagonal effective circle action 
	which on each factor is rotation around the north-south axis. It has 
	four fixed pints, and the linearization of the action induces one 
	orientation at two of the fixed points, and the other orientation at 
	the remaining two. The induced action on the boundary of an $S^{1}$-invariant 
	small ball around each of the fixed points is the Hopf action on $S^{3}$. 
	By taking equivariant connected sums at fixed points, pairing fixed points at 
	which the linearizations give opposite orientations, we obtain effective circle 
	actions with $2+2k$ fixed points on the connected sum $k(S^{2}\times S^{2})$ 
	for every $k\geq 1$. Now we have $1+k$ fixed points at which the 
	linearization induces one orientation, and $1+k$ at which it induces 
	the other orientation. Then making equivariant self-connected sums at pairs 
	of fixed points with linearizations inducing opposite orientations we 
	finally obtain a free circle action on 
	$X_{k}=k(S^{2}\times S^{2})\# (1+k)(S^{1}\times S^{3})$.

	If $k$ is odd and large enough, then there are 
	symplectic manifolds $Z_{k}$ homeomorphic (but not 
	diffeomorphic) to $k(S^{2}\times S^{2})$, see for example~\cite{HKW}. 
	By the construction given in~\cite{HKW}, we may assume that $Z_{k}$ 
	contains the Gompf nucleus of an elliptic surface. By performing 
	logarithmic transformations inside this nucleus, we can vary the 
	smooth structures on the $Z_{k}$ in such a way that the number of 
	Seiberg--Witten basic classes with numerical Seiberg--Witten 
	invariant $=\pm 1$ becomes arbitrarily large, cf.~Theorem~8.7 
	of~\cite{FSrat} and Example~3.5 of~\cite{monopoles}.
	
	Consider $Y_{k}=Z_{k}\# (1+k)(S^{1}\times S^{3})$. This is clearly 
	homeomorphic to $X_{k}$. Although the numerical Seiberg--Witten invariants 
	of $Y_{k}$ must vanish, cf.~\cite{KMT,Bourbaki}, we claim that each 
	of the basic classes with numerical Seiberg--Witten invariant $=\pm 1$ 
	on $Z_{k}$ gives rise to a monopole class on $Y_{k}$, that is the 
	characteristic class of a $\Spc$-structure for which the monopole 
	equations have a solution for every Riemannian metric on $Y_{k}$. 
	There are two ways to see this. One can extract our claim from the 
	connected sum formula~\cite{B} for the stable cohomotopy refinement of 
	Seiberg--Witten invariants introduced by Bauer and Furuta~\cite{BF}, 
	cf.~\cite{F}. Alternatively, one uses the invariant defined by 
	the homology class of the moduli space of solutions to the monopole 
	equations, as in~\cite{Bourbaki}. This means that the first homology 
	of the manifold is used, and here this is enough to obtain a non-vanishing 
	invariant. Using this invariant, our claim follows from Proposition~2.2 
	of Ozsv\'ath--Szab\'o~\cite{OS}.

	As $Y_{k}$ has non-torsion monopole classes $c$ with 
	$c^{2}=2\chi(Z_{k})+3\sigma(Z_{k})=4+4k>0$, the bound
	$$
	c^{2}\leq \frac{1}{32\pi^2}\int_{Y_{k}} s_g^{2}\dvol_{g} \ ,
	$$
	where $s_{g}$ is the scalar curvature of any Riemannian metric $g$, 
	shows that $Y_{k}$ cannot collapse with bounded scalar curvature, 
	cf.~\cite{monopoles}. {\it A fortiori} it cannot collapse with bounded 
	sectional curvature, and so its minimal volume is strictly positive. 
	
	The monopole classes we constructed on $Y_{k}$ are all generic 
	monopole classes in the sense of~\cite{monopoles}. By Lemma~2.4 of 
	{\it loc.~cit.}~each manifold has at most finitely many such classes. 
	As we can change the smooth structure to make the number of generic 
	monopole classes arbitrarily large, we have infinitely many distinct 
	smooth structures we can choose for $Y_{k}$.
	\end{proof}
The above proof also gives the following:
\begin{cor}
There are pairs of homeomorphic closed manifolds such that one collapses 
with bounded sectional curvature and the other one cannot collapse even
with bounded scalar curvature.
\end{cor}

\begin{rem}
That connected sums of manifolds with vanishing minimal volumes may have non-vanishing
minimal volumes is immediate by looking at connected sums of tori.
The manifolds $X_k$ discussed above have the property that their minimal volumes
vanish, although they are connected sums of manifolds with non-vanishing minimal volumes.
Thus the minimal volume, and even its (non-)vanishing, does not behave in a straightforward
manner under connected sums.
\end{rem}
This remark was motivated by the recent paper~\cite{PPnew} of Paternain and Petean.
After Theorem~\ref{t:2} appeared on the arXiv in~\cite{first}, 
these authors remarked on the complicated behaviour of the minimal volume under connected sums
based on some $6$-dimensional examples, see Remark~3.1 in~\cite{PPnew}.


\bigskip

\bibliographystyle{amsplain}

\begin{thebibliography}{10}

\bibitem{B1}
I.~K.~Babenko,
{\em Asymptotic invariants of smooth manifolds},
Izv.~Ross.~Akad.~Nauk.~Ser.~Mat.~{\bf 56} (1992), 707--751 (Russian);
Engl.~transl.~in Russian Acad.~Sci.~Izv.~Math.~{\bf 41} (1993), 1--38.

\bibitem{B2}
I.~K.~Babenko,
{\em Asymptotic volumes and simply connected surgeries of smooth manifolds},
Izv.~Ross.~Akad.~Nauk. Ser.~Mat.~{\bf 58} (1994), 218--221 (Russian);
Engl.~transl.~in Russian Acad.~Sci.~Izv.~Math.~{\bf 44} (1995), 427--430.

\bibitem{B3}
I.~K.~Babenko,
{\em Extremal problems of geometry, surgery on manifolds, and problems 
in group theory}, Izv.~Ross. Akad.~Nauk.~Ser.~Mat.~{\bf 59} (1995), 97--108 (Russian); 
Engl.~translation in Izv.~Math.~{\bf 59} (1995), 321--332.

\bibitem{B}
S.~Bauer, {\em A stable cohomotopy refinement of 
Seiberg--Witten invariants: II}, Invent.~math.~{\bf 155} (2004), 
21--40.

\bibitem{BF}
S.~Bauer and M.~Furuta, {\em A stable cohomotopy refinement of 
Seiberg--Witten invariants: I}, Invent.~math.~{\bf 155} (2004), 1--19.

\bibitem{besse}
A.~L.~Besse, {\sl Einstein Manifolds}, Springer Verlag 1987.

\bibitem{Be}
L.~Bessi\`eres, {\em Un th\'eor\`eme de rigidit\'e diff\'erentielle}, 
Comment.~Math.~Helv.~{\bf 73} (1998), 443--479.

\bibitem{BCG0}
G.~Besson, G.~Courtois et S.~Gallot, {\em Volume et entropie minimale 
des espaces localement sym\'etriques}, Invent.~math.~{\bf 103}
(1991), 417--445.

\bibitem{BCG}
G.~Besson, G.~Courtois et S.~Gallot, {\em Entropies et 
rigidit\'es des espaces localement sym\'etriques de courbure
strictement n\'egative}, Geom.~Func.~Analysis (GAFA) {\bf 5}
(1995), 731--799.

\bibitem{BC}
R.~L.~Bishop and R.~J.~Crittenden, 
{\sl Geometry of Manifolds}, Academic Press, New York 1964; reprinted 
by AMS Chelsea Publishing, Providence 2001.


\bibitem{br}
M.~Brunnbauer, {\em Homological invariance for asymptotic invariants and systolic inequalities}, 
Geom.~Func.~Analysis (GAFA) {\bf 18} (2008), 1087--1117.

\bibitem{BIS}
M.~Brunnbauer, M.~Ishida and P.~Su\'arez-Serrato, {\em An essential relation between
Einstein metrics, volume entropy, and exotic smooth structures}, 
Math.~Research Letters {\bf 16} (2009), 503--514.

\bibitem{Cheeger}
J.~Cheeger, {\em A lower bound for the smallest eigenvalue of the Laplacian}, 
in {\sl Problems in analysis} (Papers dedicated to Salomon Bochner, 1969), pp. 195--199. 
Princeton Univ. Press, Princeton, N. J., 1970.

\bibitem{FSrat}
R.~Fintushel and R.~J.~Stern, {\em Rational blowdowns of smooth $4$--manifolds}, 
J.~Differential Geometry {\bf 46} (1997), 181--235.

\bibitem{F}
M.~Furuta, {\it Private communication}, 2003.

\bibitem{gromov}
M.~Gromov, {\em Volume and bounded cohomology}, 
Publ.~Math.~I.H.E.S. {\bf 56} (1982), 5--99.

\bibitem{HKW}
B.~Hanke, D.~Kotschick and J.~Wehrheim,
{\em Dissolving four-manifolds and positive scalar curvature},
Math.~Zeit.~{\bf 245} (2003), 545--555.

\bibitem{Hil}
J.~A.~Hillman, {\em Flat $4$-manifold groups}, New Zealand 
J.~Math.~{\bf 24} (1995), 29--40.

\bibitem{HT}
N.~J.~Hitchin, {\em Compact four--dimensional Einstein manifolds}, 
J.~Differential Geometry {\bf 9} (1974), 435--441.

\bibitem{KKM}
J.~Kedra, D.~Kotschick and S.~Morita,
{\em Crossed flux homomorphisms and vanishing theorems for flux groups}, 
Geom.~Func.~Analysis (GAFA) {\bf 16} (2006), 1246--1273.

\bibitem{Bourbaki}
D.~Kotschick, {\em The Seiberg-Witten invariants of symplectic 
four--manifolds}, S\'eminaire Bourbaki, 48\`eme ann\'ee, 1995-96, 
no.~812, Ast\'erisque {\bf 241} (1997), 195--220.

\bibitem{K}
D.~Kotschick, {\em Einstein metrics and smooth structures},
Geometry \& Topology {\bf 2} (1998), 1--10.

\bibitem{GHT}
D.~Kotschick, {\em On the Gromov--Hitchin--Thorpe inequality},
C.~R.~Acad.~Sci.~Paris {\bf 326} (1998), 727--731.

\bibitem{monopoles}
D.~Kotschick,
{\em Monopole classes and Einstein metrics},
Intern.~Math.~Res.~Notices {\bf 2004} no.~12 (2004), 593--609.

\bibitem{first}
D.~Kotschick, {\em Entropies, volumes, and Einstein metrics}, 
Preprint arXiv:math/0410215 v1 [math.DG] 8 Oct 2004.

\bibitem{KMT}
D.~Kotschick, J.~W.~Morgan and C.~H.~Taubes, {\em Four-manifolds 
without symplectic structures but with non-trivial Seiberg--Witten 
invariants}, Math. Research Letters {\bf 2} (1995), 119--124.

\bibitem{LW}
F.~Ledrappier and X.~Wang, {\em An integral formula for the volume 
entropy with applications to rigidity},  J.~Differential Geometry {\bf 85} (2010), 461--477.

\bibitem{M}
A.~Manning,
{\em Topological entropy for geodesic flows}, 
Ann.~of Math.~{\bf 110} (1979), 567--573.

\bibitem{M2}
A.~Manning,
{\em More topological entropy for geodesic flows}, 
In: Dynamical Systems and Turbulence, Springer LNM {\bf 898} (1981), 
234--249.

\bibitem{Mi}
J.~Milnor, {\em A note on curvature and the fundamental group}, 
J.~Differential Geometry {\bf 2} (1968), 1--7.

\bibitem{OS}
P.~Ozsv\'ath and Z.~Szab\'o, {\em Higher type adjunction inequalities 
in Seiberg--Witten theory}, J.~Differential Geometry {\bf 55} (2000), 385--440.
 
\bibitem{P}
G.~P.~Paternain, 
{\sl Geodesic Flows},
Progress in Mathematics vol.~{\bf 180}, Birk\"auser Verlag 1999.

\bibitem{PPMRL}
G.~P.~Paternain and J.~Petean,
{\em Einstein manifolds of non-negative sectional curvature and entropy},
Math.~Research Letters~{\bf 7} (2000), 503--515.

\bibitem{PP}
G.~P.~Paternain and J.~Petean,
{\em Minimal entropy and collapsing with curvature bounded from 
below},
Invent.~math. {\bf 151} (2003), 415--450.

\bibitem{PPnew}
G.~P.~Paternain and J.~Petean,
{\em Collapsing manifolds obtained by Kummer-type constructions},
Trans.~Amer.~Math. Soc.~{\bf 361} (2009), 4077--4090.

\bibitem{S}
A.~Sambusetti, {\em An obstruction to the existence of Einstein
metrics on $4$--manifolds}, Math.~Ann.~{\bf 311} (1998), 533--547.

\bibitem{Sch}
R.~Schmidt, {\sl Spectral geometry, asymptotic invariants, and geometric 
group theory}, Diplomarbeit M\"unchen 2009.

\bibitem{wagner}
M.~Wagner, {\sl \"Uber die Klassifikation flacher Riemannscher 
Mannigfaltigkeiten}, Diplomarbeit Basel 1997.

\end{thebibliography}

\bigskip

\end{document}